\documentclass{article}
\catcode`\@=11 \@addtoreset{equation}{section}

\catcode`\@=12

\newcommand{\mathbb}{\mathbf}

\newtheorem{Theorem}{Theorem}[section]
\newtheorem{Proposition}{Proposition}[section]
\newtheorem{Lemma}{Lemma}[section]
\newtheorem{Corollary}{Corollary}[section]
\newtheorem{Definition}{Definition}[section]
\newtheorem{Example}{Example}[section]

\newcommand{\bTheorem}[1]{
\begin{Theorem} \label{T#1} }
\newcommand{\eT}{\end{Theorem}}

\newcommand{\bProposition}[1]{
\begin{Proposition} \label{P#1}}
\newcommand{\eP}{\end{Proposition}}

\newcommand{\bLemma}[1]{
\begin{Lemma} \label{L#1} }
\newcommand{\eL}{\end{Lemma}}

\newcommand{\bCorollary}[1]{
\begin{Corollary} \label{C#1} }
\newcommand{\eC}{\end{Corollary}}

\newcommand{\bFormula}[1]{
\begin{equation} \label{#1}}
\newcommand{\eF}{\end{equation}}

\newcommand{\Ov}[1]{\overline{#1}}
\newcommand{\DC}{C^\infty_c}

\newcommand{\vr}{\varrho}
\newcommand{\vre}{\vr_n}

\newcommand{\vue}{\vu_n}

\newcommand{\vu}{\vc{u}}
\newcommand{\vc}[1]{{\bf #1}}

\newcommand{\Div}{{\rm div}_x}
\newcommand{\Grad}{\nabla_x}

\newcommand{\tn}[1]{\mbox {\F #1}}
\newcommand{\dx}{{\rm d} {x}}
\newcommand{\dt}{{\rm d} t }

\newcommand{\dxdt}{\dx \ \dt}
\newcommand{\intO}[1]{\int_{\Omega} #1 \ \dx}

\newcommand{\ep}{\varepsilon}

\font\F=msbm10 scaled 1000

\date{}
\begin{document}


\title{Compressible fluid flows driven by stochastic forcing}
\author{Eduard Feireisl \thanks{The work was supported by Grant 201/09/
0917 of GA \v CR as a part of the general research programme of the
Academy of Sciences of the Czech Republic, Institutional Research
Plan RVO: 67985840.}\and Bohdan Maslowski \thanks{The work was supported by Grant P201/10/0752 of  GA \v CR.} \and Anton{\' \i}n Novotn\'
y \thanks{The work was supported the research programme of the
Academy of Sciences of the Czech Republic, Institutional Research
Plan AV0Z10190503.}}

\maketitle

\bigskip

\centerline{Institute of Mathematics of the Academy of Sciences of
the Czech Republic} \centerline{\v Zitn\' a 25, 115 67 Praha 1,
Czech Republic}

\medskip

\centerline{Faculty of Mathematics and Physics, Charles University in Prague}
\centerline{Sokolovska 83, Praha 8, Czech Republic}

\medskip

\centerline{IMATH, Universit\' e du Sud Toulon-Var, BP 132, 839 57
La Garde, France}

\medskip

\begin{abstract}

We consider the Navier-Stokes system describing the motion of a
compressible barotropic fluid driven by stochastic external forces.
Our approach is semi-deterministic, based on solving the system for
each fixed representative of the random variable and applying an
abstract result on measurability of multi-valued maps. The class of
admissible driving forces includes the (temporal) white noise and
the random kicks, considered recently in the context of
incompressible fluid models.

\end{abstract}


\bigskip

\section{Problem formulation}
\label{p}

We consider the Navier-Stokes system governing the time evolution of
the density $\vr$ and the velocity $\vu$ of a compressible viscous
fluid, driven by a stochastic external force that can be formally written in the form:
\bFormula{p1}
{\rm d} \vr + \Div (\vr \vu ) \ {\rm d}t  = 0,
\eF
\bFormula{p2}
{\rm d} (\vr \vu) + \Big( \Div (\vr \vu \otimes \vu) + \Grad p(\vr)
- \Div \tn{S} \Big)\ {\rm d}t   = \vr {\rm d} \vc{w}.
\eF
The quantities $\vr = \vr(t,x,\omega)$, $\vu = \vu(t,x,\omega)$ are
functions of the time $t \in (0,T)$, the spatial position $x \in
\Omega$, and $\omega \in {\mathcal O} = \{ {\mathcal O}, \mathcal B
, \mu \}$, where ${\mathcal O}$ is a topological probability space,
with the family of Borel sets ${\mathcal B}$, and a regular
probability measure $\mu$. System (\ref{p1}), (\ref{p2}) is
supplemented with the standard no-slip boundary condition
\bFormula{p3}
\vu|_{\partial \Omega} = 0,
\eF
and the initial conditions
\bFormula{p4}
\vr(0, \cdot) = \vr_0 , \ \vr \vu (0,\cdot) = (\vr \vu)_0.
\eF

The symbol $\tn{S}$ denotes the viscous stress determined by
Newton's rheological law
\bFormula{p3a}
\tn{S} = \nu \Big( \Grad \vu + \Grad^t \vu - \frac{2}{3} \Div \vu
\Big) + \eta \Div \vu \tn{I},\ \nu > 0, \ \eta \geq 0,
\eF
$p = p(\vr)$ is the pressure, and the perturbation $\vc{w}$ is a
random variable represented for a.a. $\omega$ by a bounded function,
sufficiently regular with respect to the spatial variable $x \in
\Omega$.

There is already a substantial amount of work on the \emph{incompressible} stochastic Navier-Stokes system,
see the surveys by Bensoussan \cite{Bens}, Mattingly \cite{Matt},
and the references cited therein. In 3D general existence results for equations with Gaussian noise have been obtained for the so-called weak solutions in probabilistic sense (or martingale solutions), cf. \cite{GM}, \cite{FL} or \cite{FR}. Uniqueness of solutions in 3D is open just as in the deterministic case, but working carefully with the so-called Markov selections (\cite{DD}, \cite{FR}) it was possible to prove existence, uniqueness, ergodicity and strong mixing of the invariant measure (stationary solution) assuming certain nondegeneracy of the noise (for a general result on Markov selections including also the 3D stochastic Navier-Stokes equation cf. \cite{GRZ}). Analogous results (existence of martingale solutions and Markov selection) for 3D Navier-Stokes equations with jumps (driven by L\'evy noise) have been recently proved by Dong and Zhai \cite{DZ}.

 Much less seems to be known for \emph{compressible} fluid flows. To the best of our knowledge, the only available results concerning the compressible stochastic Navier-Stokes system concern the 1-D case,
see Tornatore , Fujita Yashima \cite{TorYasB}, \cite{TorYasA}, and
the rather special peridic 2-D case examined by Tornatore \cite{Torn} by means of
the existence theory developed by Vaigant and Kazhikhov \cite{VAKA}.

Similarly to the seminal paper by Bensoussan and Temam \cite{BenTem}, our approach is semi-deterministic based on the concept of weak solutions to the Navier-Stokes system, in the framework of the existence theory developed by Lions \cite{LI4}, with the extension of \cite{FNP}. The class of the weak solutions is sufficiently robust to incorporate the driving forces
$\vc{f}$ of very low regularity with respect to the time, in particular, we may consider
\[
\vc{f} = \vc{f}(t,x) = \partial_t \vc{w}(t,x), \mbox{as soon as} \ \vc{w} \ \mbox{is Lipschitz with respect to the spatial
variable}\ x \in \Omega.
\]
Thus the stochastic problem may be solved ``pathwise'' for any individual choice of the random element $\omega$. It is also shown that under some mild restrictions on the paths of driving process ( the time integral of the noise) the solution is a "stochastic process", i.e. it is a (measurable) random variable with values in the paths space. As an example we consider the L\'evy noise , which is important because it represents a general stationary incorrelated noise driving the paths which are stochastically continuous (cf. Section 5 for details).

\subsection{Weak formulation}
\label{wf}

Taking advantage of the specific form of (\ref{p1}), (\ref{p2}), we
may formally write
\[
{\rm d} (\vr \vu) - \vr {\rm d} \vc{w} = {\rm d} [ \vr ( \vu -
\vc{w})] + {\rm d} \vr \vc{w} = {\rm d} [ \vr ( \vu - \vc{w})] -
\Div (\vr \vu) \vc{w} \ {\rm d}t.
\]
Accordingly, the product $\Div (\vr \vu)$ may be interpreted in the
sense of distributions provided $\vc{w}$ ranges in the Sobolev space
$W^{1,q}_0(\Omega;R^3)$ for a sufficiently large exponent $q$.

We say that a pair of random variables $\vr$, $\vu$,
\bFormula{random1}
\vr = \vr(t,x, \omega), \ \vr(\cdot, \omega) \in C([0,T];L^1(\Omega)),
\eF
\bFormula{random2}
\vu = \vu(t,x,\omega),\
\vu(\cdot, \omega) \in L^2(0,T; W^{1,2}_0(\Omega;R^3)) \ \mbox{for a.a.}\ \omega \in \mathcal{O},
\eF
is a \emph{weak solution} to the problem (\ref{p1} -
\ref{p4}), if the following relations hold:

\begin{enumerate}
\item

\centerline{\textsc{Equation of continuity:}}

\bFormula{ws1}
\vr \geq 0,\  \vr \in L^\infty(0,T:L^\gamma(\Omega)) \ \mbox{for a certain}\ \gamma > 1
\ \mbox{for a.a.}\ \omega \in
{\mathcal O},
\eF
and the family of integral identities
\bFormula{ws2}
\intO{ \Big( \vr (\tau, x, \omega)  + b(\vr(\tau, x, \omega) ) \Big) \varphi(x)  } - \intO{ \Big( \vr_0 (x, \omega)+ b(\vr_0(x, \omega)) \Big)
\varphi(x) }
\eF
\[
= \int_0^\tau \intO{ \Big( \vr (t,x,\omega) \vu  (t,x,\omega) \cdot
\Grad \varphi(x)   + \Big( b(\vr(t,x,\omega) ) - b'(\vr(t,x,\omega)) \vr(t,x,\omega) \Big) \Div \vu (t,x,\omega) \varphi (x) \Big) } \dt
\]
is satisfied
for any test function $\varphi \in C^1(\Ov{\Omega})$, any $b \in \DC[0,\infty)$, and a.a.
$\omega \in \mathcal{O}$;
\item

\centerline{\textsc{Momentum equation:}}

\bFormula{ws3}
\vr \vu \in L^\infty(0,T: L^{2\gamma/ (\gamma + 1)}(\Omega;R^3)),\
\vr (\vu - \vc{w}) \in C_{\rm weak}([0,T] ; L^{2\gamma /(\gamma +
1)}(\Omega;R^3)) \ \mbox{for a.a}\ \omega \in {\mathcal O},
\eF
\bFormula{ws4}
p(\vr) \in L^\infty(0,T; L^1(\Omega) \cap L^q((0,T) \times \Omega) \ \mbox{for some} \ \mbox{for a.a.}\ \omega \in
\mathcal{O},
\eF
and the integral identity
\bFormula{ws5}
\intO{ \vr(\tau, x, \omega) \Big( \vu(\tau, x, \omega) - \vc{w}(\tau, x, \omega) \Big) \cdot \varphi(x) } -
\intO{ (\vr \vu)_0 (x, \omega) \cdot \varphi(x) }
\eF
\[
= \int_0^\tau \intO{ \Big( \vr(t,x, \omega) \vu (t,x, \omega) \otimes \vu
(t,x, \omega) : \Grad \varphi(x) + p(\vr(t,x, \omega)) \Div \varphi(x) \Big) }
\ \dt
\]
\[
 - \int_0^\tau \intO{ \Big(
\tn{S}(\Grad \vu (t,x,\omega)) : \Grad \varphi(x)  + \vr
(t,x,\omega)\vu (t,x,\omega) \cdot \Grad \Big( \vc{w}(t,x, \omega)
\varphi(x) \Big) \Big)} \ \dt
\]
holds
for all $\varphi \in \DC(\Omega;R^3)$, and for a.a. $\omega \in
\mathcal{O}$.

\item

\centerline{\textsc{Energy inequality}}

\bFormula{ws6}
\intO{ \left( \frac{1}{2} \vr(\tau,x,\omega) | \vu(\tau,x,\omega) -
\vc{w}(\tau, x, \omega)  |^2 + P(\vr (\tau, x, \omega ) ) \right) }
\eF
\[
+ \int_s^\tau \intO{ \tn{S} (\Grad \vu (t,x,\omega)): \Grad \vu (t,x,\omega)
}
\]
\[
\leq \intO{ \left( \frac{1}{2} \vr(s,x,\omega) | \vu(s,x,\omega) -
\vc{w}(s, x, \omega)  |^2 + P(\vr (s, x, \omega ) ) \right) }
\]
\[
+\int_s^\tau \intO{ \Big( \tn{S}(\Grad \vu) : \Grad \vc{w} - \vr \vu
\otimes \vu : \Grad \vc{w} - p(\vr) \Div \vc{w} + \frac{1}{2} \vr
\vu \cdot \Grad |\vc{w} |^2 \Big)(t,x,\omega) } \ \dt
\]
for a.a. $\tau \geq  s \geq 0$ including $s = 0$, and a.a. $\omega \in \mathcal{O}$,
with
\bFormula{ws7}
P(\vr) = \vr \int_1^\vr \frac{p(z)}{z^2} \ {\rm d}z.
\eF
\end{enumerate}

{\bf Remark \ref{wf}.1} {\it The family of integral identities
(\ref{ws2}) represents a weak formulation of the renormalized
equation of continuity introduced by DiPerna and Lions \cite{DL}.
The more standard, but in fact equivalent, ``distributional''
formulation of (\ref{ws2}), (\ref{ws5}), (\ref{ws6}) will be given
in Section \ref{s} below. }

\medskip

We can check easily that for the energy inequality (\ref{ws6}) to
yield uniform bounds on the family of solutions, it is necessary
that
\bFormula{ws7a}
{\rm ess} \sup_{t \in (0,T)} \| \vc{w} (t, \cdot, \omega) \|_{W^{1,\infty}_0(\Omega;R^3)} \leq c(\omega) \ \mbox{for a.a.}\ \omega \in \mathcal{O},
\eF
in particular, certain regularity with respect to the spatial variable is
needed.

\subsection{Main results}

Our main goal is to develop an existence theory for the Navier-Stokes system
(\ref{ws1} - \ref{ws6}), where the solutions $\vr(\cdot, \omega)$, $\vu(\cdot, \omega)$ are random variables ranging in the function spaces specified through (\ref{random1}), (\ref{random2}). Analogously to
Bensoussan and Temam \cite{BenTem}, our approach is based on:

\begin{itemize}
\item a general \emph{stability} result for solutions of (\ref{ws1}
- \ref{ws6}), where the forcing term satisfies only (\ref{ws7a});
\item application of the abstract \emph{measurability theorem} for multivalued maps,
similar to Bensoussan and Temam \cite[Theorem 3.1]{BenTem}.
\end{itemize}

By \emph{stability} we mean that a sequence of weak solutions to
(\ref{ws1} - \ref{ws6}), with pre-compact data $\vr_0$, $(\vr
\vu)_0$, $\vc{w}$, admits a subsequence that converges (weakly) to
another solution of the same problem. Besides the nowadays standard
ingredients of the existence theory (see \cite{FNP}, Lions
\cite{LI4}), the proof is based on compactness of certain
commutators in the spirit of Coifman and Meyer \cite{COME}.

The paper is organized as follows. In Section \ref{s}, we recall the standard definition of weak solutions to the Navier-Stokes system and prove a general stability results assuming low regularity of the driving force represented by $\vc{w}$. This result is exploited in Section \ref{e}, where we show existence of the weak solutions for system (\ref{ws1} - \ref{ws6}), with a fixed (irregular) force $\vc{w}$. Section \ref{ps} is devoted to the analysis of the associated stochastic system. We introduce the suitable function-spaces framework and use the abstract measurability theorem to prove that the weak solution generate a random variable. Several concrete examples
of the driving force $\vc{w}$ are discussed in Section \ref{sk}.

\section{Weak sequential stability for irregular forcing}
\label{s}

In order to fix ideas, we suppose that the pressure $p$ belongs to the class
$C[0,\infty) \cap C^2(0,\infty)$ and satisfies
\bFormula{w4}
p(0) = 0, \ p'(\vr) > 0 \ \mbox{for}\ \vr > 0, \ \lim_{\vr \to
\infty} \frac{p'(\vr)}{\vr^{\gamma - 1}} = p_\infty > 0,\ \gamma >
3/2.
\eF

In this section, we drop the parameter $\omega$ and fix a function $\vc{w}$.
Accordingly, we may rewrite the weak formulation of the problems as follows:

\begin{enumerate}
\item

\centerline{\textsc{Equation of continuity (renormalized):}}

\bFormula{Fws2}
\int_0^T \intO{ \left( \Big(\vr + b(\vr) \Big) \partial_t \varphi + \Big( \vr + b(\vr) \Big) \vu \cdot \Grad \varphi
\right) } \ \dt
\eF
\[
= \int_0^T \intO{  \Big( b'(\vr)\vr - b(\vr) \Big) \Div \vu \Big) \varphi } \ \dt - \intO{ \Big( \vr_0 + b(\vr_0) \Big) \varphi (0, \cdot) }
\]
for any test function $\varphi \in \DC([0,T) \times \Ov{\Omega} )$, and any $b \in \DC[0,\infty)$;
\item

\centerline{\textsc{Momentum equation:}}

\bFormula{Fws5}
\int_0^T \intO{ \Big( \vr (\vu - \vc{w})  \cdot \partial_t \varphi + \vr \vu \otimes
\vu : \Grad \varphi + p(\vr) \Div \varphi \Big)  } \ \dt
\eF
\[
\int_0^T \intO{ \Big( \tn{S}(\Grad \vu ) : \Grad \varphi(x)  + \vr
\vu \cdot \Grad ( \vc{w} \cdot \varphi ) \Big)} \ \dt - \intO{ (\vr
\vu)_0 \cdot \varphi (0,\cdot) }
\]
for all $\varphi \in \DC([0,T) \times \Omega;R^3)$;

\item

\centerline{\textsc{Energy inequality:}}

\bFormula{Fws6}
- \int_0^T \intO{ \left( \frac{1}{2} \vr | \vu -
\vc{w}  |^2 + P(\vr  ) \right)} \ \partial_t \psi \ \dt
+ \int_0^T \intO{ \tn{S} (\Grad \vu ): \Grad \vu
} \ \psi \ \dt
\eF
\[
\leq \psi(0) \intO{ \left( \frac{1}{2} \frac{ |(\vr \vu)_0 |^2 }{\vr_0}  + P(\vr_0)  \right) }
\]
\[
+\int_0^T \intO{ \Big( \tn{S}(\Grad \vu) : \Grad \vc{w} - \vr \vu
\otimes \vu : \Grad \vc{w} - p(\vr) \Div \vc{w} + \frac{1}{2} \vr
\vu \cdot \Grad |\vc{w} |^2 \Big) } \ \psi \ \dt
\]
for any $\psi \in \DC[0, T)$, $\psi \geq 0$.
\end{enumerate}

The heart of the paper is the following result that may be of independent interest.

\bProposition{s1}
Let $\Omega \subset R^3$ be a bounded Lipschitz domain. Suppose that the
pressure $p$ satisfies (\ref{w4}). Let $\{ \vc{w}_n \}_{n=1}^\infty$ be a sequence of functions,
\bFormula{s3}
{\rm ess} \sup_{t \in (0,T)} \| \vc{w}_n (t, \cdot) \|_{W^{1,\infty}_0(\Omega;R^3)} \leq C_w,
\vc{w}_n \to \vc{w} \ \mbox{in}\ L^1(0,T; W^{1,1}(\Omega; R^3)).
\eF
Let $\{ \vr_n, \vu_n \}_{n=1}^\infty$ be a sequence of weak solutions of the Navier-Stokes system in $(0,T) \times \Omega$, driven by $\vc{w}_n$, and emanating from the initial data $\vr_{0,n}$, $(\vr \vu)_{0,n}$ such that
\bFormula{s1}
\vr_{0,n} \geq 0,\ \intO{ \vr_{0, n} } = M > 0, \ \| \vr_{0,n} \|_{L^\gamma(\Omega)} \leq E_0 ,\ \intO{ \frac{ |(\vr
\vu)_{0,n}|^2 }{\vr_{0,\ep} } }  \leq E_0,
\eF
\bFormula{s2}
\vr_{0, n} \to \vr_0 \ \mbox{in} \ L^\gamma(\Omega),\ (\vr
\vu)_{0, n} \to (\vr \vu)_0 \ \mbox{weakly in}\ L^1(\Omega), \
\intO{ \frac{ |(\vr \vu)_{0,n}|^2 }{\vr_{0,n} } } \to \intO{
\frac{ |(\vr \vu)_{0}|^2 }{\vr_{0} } }.
\eF

Then, at least for suitable subsequences,
\bFormula{s4}
\vr_n \to \vr \ \mbox{in}\ C_{\rm weak}([0,T]; L^\gamma(\Omega))
\ \mbox{and in} \ L^1((0,T) \times \Omega),
\eF
\bFormula{s5}
\vu_n \to \vu \ \mbox{weakly in} \ L^2(0,T; W^{1,2}_0(\Omega;R^3)),
\eF
\bFormula{s6}
\vr_n(\vu_n - \vc{w}_n ) \to \vr (\vu - \vc{w}) \ \mbox{in}\
C_{\rm weak}([0,T]; L^{2\gamma/(\gamma + 1)}(\Omega; R^3)),
\eF
where $\vr$, $\vu$ is a weak solution of the same problem, with the driving force $\vc{w}$,
and the initial data $\vr_0$, $(\vr \vu)_0$.
\eP

{\bf Remark \ref{s}.1} {\it In hypothesis (\ref{s1}), we tacitly assume that
$(\vr \vu)_{0,n} (x) = 0$ on the set where $\vr_{0,n} (x) = 0$.
}

\medskip

The remaining part of this section is devoted to the proof of Proposition \ref{Ps1}.

\subsection{Uniform bounds}

Since the total mass $M$ is a constant of motion, we get
\bFormula{ub1}
 \| \vre(t, \cdot) \|_{L^1(\Omega)} = M \ \mbox{for any}\ t \in [0,T].
\eF
Consequently, in accordance with hypothesis (\ref{s3}),
\[
{\rm ess} \sup_{t \in (0,T)} \intO{ \vre |\vc{w}_n|^2 } \leq c(M, C_w);
\]
whence the standard Gronwall argument can be used in the energy inequality
(\ref{Fws6}) to obtain
\bFormula{ub2}
{\rm ess} \sup_{t \in (0,T)} \| \sqrt{\vre} \vue (t, \cdot) \|_{L^2(\Omega;R^3)} \leq c(
M, C_w, E_0,T),
\eF
\bFormula{ub3}
{\rm ess} \sup_{t \in (0,T)} \| \vre (t, \cdot) \|_{L^\gamma(\Omega)} \leq
c(M, C_w, E_0, T),
\eF
and, by virtue of Korn's inequality,
\bFormula{ub4}
\int_0^T \| \vue (t, \cdot) \|^2_{W^{1,2}_0(\Omega;R^3)} \leq c(M,C_w, E_0,T).
\eF

The estimates (\ref{ub1} - \ref{ub4}) are exactly the same as the {\it a priori} bounds available for the compressible Navier-Stokes system, cf. \cite{FNP}. In particular, we may infer, at least for suitable subsequences,
that
\bFormula{ub5}
\vre \to \vr \ \mbox{in}\ C_{\rm weak}([0,T]; L^\gamma (\Omega)),
\eF
\bFormula{ub6}
\vue \to \vu \ \mbox{weakly in} \ L^2(0,T; W^{1,2}_0(\Omega;R^3)),
\eF
and
\bFormula{ub8}
\vre \vue \to \vr \vu \ \mbox{weakly-(*) in}\ L^\infty(0,T;
L^{2\gamma/(\gamma + 1)}(\Omega; R^3).
\eF
Moreover, a short inspection of the momentum equation (\ref{Fws5}) yields
\bFormula{ub9}
\vre (\vue - \vc{w}_n ) \to \vr (\vu - \vc{w}) \ \mbox{in}\
C_{\rm weak} ([0,T]; L^{2 \gamma/ (\gamma + 1)}(\Omega;R^3)).
\eF

Finally, combining (\ref{ub8}), (\ref{ub9}) with hypothesis (\ref{s3}), we obtain that
\bFormula{ub10}
(\vre \vue) (t, \cdot) \to (\vr \vu) (t, \cdot) \ \mbox{weakly in}\
L^{2\gamma / (\gamma + 1)} (\Omega;R^3) \ \mbox{for a.a.}\ t \in (0,T);
\eF
in particular, as $\frac{2 \gamma}{\gamma + 1} > \frac{6}{5}$,
\bFormula{ub11}
\vre \vue \to \vr \vu \ \mbox{in}\ L^q(0,T; W^{-1,2}(\Omega;R^3))
\ \mbox{for any}\ 1 \leq q < \infty.
\eF
Relations (\ref{ub6}), (\ref{ub11}), together with the bounds established in
(\ref{ub2}), (\ref{ub3}), give rise to
\bFormula{ub12}
\vre \vue \otimes \vue \to \vr \vu \otimes \vu \ \mbox{weakly in}
\ L^2(0,T; L^{6 \gamma /(4 \gamma + 3)}(\Omega;R^{3 \times 3})).
\eF

\subsection{Pressure estimates}

In view of the uniform bounds established in (\ref{ub1} - \ref{ub4}), the pressure $p$ can be estimated exactly as in \cite{FP13}, specifically, by means of the quantities
\[
\varphi = \psi (t) {\cal B} \left[ b(\vre) - \frac{1}{|\Omega|} \intO{
b(\vre) } \right], \ \psi \in \DC(0,T), \ b(\vr) \approx \vr^\nu
\]
used as test functions in the momentum equation (\ref{Fws5}),
where ${\cal B} \approx {\rm div}^{-1}$ is the so-called \emph{Bogovskii operator}, see Bogovskii \cite{BOG}, Galdi \cite{GAL}. Indeed the extra terms that appear in (\ref{Fws5}) in comparison with the situation treated in
\cite{FP13} are regular as a consequence of hypothesis (\ref{s3}). Consequently, exactly as in \cite{FP13}, we deduce that
\bFormula{ub13}
\int_0^T \intO{ p(\vre) \vre^\beta } \leq c(M,C_w, E_0,T)
\ \mbox{for a certain} \ \beta = \beta(M, C_w, E_0) > 0.
\eF

\subsection{Pointwise convergence of the densities}

In order to complete the proof of Proposition \ref{Ps1},
we establish the following crucial result:
\bFormula{c1}
\vre \to \vr \ \mbox{a.a. in} \ (0,T) \times \Omega.
\eF
Obviously, the a.a. pointwise convergence of the densities is a key ingredient of the existence theory for the Navier-Stokes system, see
\cite{FNP}, Lions \cite{LI4}. Here, the situation is even more delicate because of the presence of the new terms containing $\vc{w}$.

Following the arguments of \cite{FNP}, we use
\[
\varphi = \psi(t) \eta(x) \Grad \Delta^{-1}[ \xi b(\vre) ],\
\psi \in \DC(0,T), \ \eta, \xi \in \DC(\Omega)
\]
as test functions in the momentum equation (\ref{Fws5}), and, similarly,
we take
\[
\varphi = \psi \eta (x) \Grad \Delta^{-1} [\xi \Ov{b(\vr)}]
\]
in the limit equation
\bFormula{c2}
\int_0^T \intO{ \Big( \vr (\vu - \vc{w} ) \cdot \partial_t \varphi +
\vr \vu \otimes \vu : \Grad \varphi + \Ov{p(\vr)} \Div \varphi \Big) } \
\dt
\eF
\[
= \int_0^T \intO{ \tn{S}(\Grad \vu) : \Grad \varphi } \ \dt +
\int_0^T \intO{ \vr \vu \cdot \Grad ( \vc{w} \cdot \varphi ) } \ \dt
- \intO{ (\vr \vu)_0 \cdot \varphi (0, \cdot) },
\]
where we have used $\Ov{h(\vr)}$ to denote a weak limit of the compositions
$h(\vre)$.

Following step by step the arguments  of \cite{FNP} we deduce that
\bFormula{c3}
\lim_{n \to \infty} \int_0^T \int_{R^3} \xi \eta \Big(
\xi p(\vre) b(\vre) - \tn{S}_\ep : \Grad \Delta^{-1} \Grad [ \xi b(\vre) ]
\Big)
\ \dxdt
\eF
\[
- \int_0^T \int_{R^3} \xi \eta \Big(
\xi \Ov{p(\vr)}\ \Ov{b(\vr)} - \tn{S} : \Grad \Delta^{-1} \Grad [ \xi \Ov{b(\vr)} ]
\Big)
\ \dxdt
\]
\[
=
\lim_{n \to \infty} \int_0^T \int_{R^3} \psi \vue \cdot \Big(
\Grad \Delta^{-1} \Grad [ \xi b(\vre)] \cdot \eta \vre \vue -
\xi b(\vre) \Grad \Delta^{-1} \Grad [\eta \vre \vue ] \Big) \ \dxdt
\]
\[
-\int_0^T \int_{R^3} \psi \vu \cdot \Big(
\Grad \Delta^{-1} \Grad [ \xi \Ov{b(\vr)}] \cdot \eta \vr \vu -
\xi \Ov{b(\vr)} \Grad \Delta^{-1} \Grad [\eta \vr \vu ] \Big) \ \dxdt
\]
\[
+ \lim_{n \to \infty} \int_0^T \int_{R^3} \left[ \psi \vre \vue
\cdot \Grad \Big( \eta \vc{w}_n \cdot \Grad \Delta^{-1} [ \xi
b(\vre) ] \Big) + \vre \vc{w}_n \cdot \partial_t \Big( \psi \eta
\Grad \Delta^{-1} [ \xi b(\vre) ] \Big) \right] \dxdt
\]
\[
- \int_0^T \int_{R^3} \left[ \psi \vr \vu \cdot \Grad \Big( \eta
\vc{w} \cdot \Grad \Delta^{-1} [ \xi \Ov{b(\vr)} ] \Big) + \vr
\vc{w} \cdot \partial_t \Big( \psi \eta \Grad \Delta^{-1} [ \xi
\Ov{b(\vr)} ] \Big) \right] \dxdt.
\]

\subsubsection{Weak stability of the effective viscous flux}

We show that the right-hand side of (\ref{c3}) vanishes, in particular, we recover the so-called weak continuity of the effective viscous flux, first established by Lions \cite{LI4}:
\bFormula{wws1}
\Ov{ p(\vr) b(\vr) } - \left( \frac{4}{3} \nu + \lambda \right)
\Ov{ b(\vr) \Div \vu } = \Ov{ p(\vr)} \ \Ov{ b(\vr) } -\left( \frac{4}{3} \nu + \lambda \right)
\Ov{ b(\vr) } \Div \vu
\eF
for any bounded function $b$.
In order to see (\ref{wws1}), observe first that we can use the same arguments
as in \cite{FNP} to show that
\bFormula{conv}
\lim_{n \to \infty} \int_0^T \int_{R^3} \psi \vue \cdot \Big(
\Grad \Delta^{-1} \Grad [ \xi b(\vre)] \cdot \eta \vre \vue -
\xi b(\vre) \Grad \Delta^{-1} \Grad [\eta \vre \vue ] \Big) \ \dxdt
\eF
\[
= \int_0^T \int_{R^3} \psi \vu \cdot \Big(
\Grad \Delta^{-1} \Grad [ \xi \Ov{b(\vr)}] \cdot \eta \vr \vu -
\xi \Ov{b(\vr)} \Grad \Delta^{-1} \Grad [\eta \vr \vu ] \Big) \ \dxdt.
\]
Indeed the only difference between \cite{FNP} and the present situation is that the convergence
\[
\vre \vue \to \vr \vu \ \mbox{in} \ C_{\rm weak}([0,T]; L^{2\gamma/(\gamma + 1)}(\Omega)),
\]
valid in \cite{FNP}, is replaced by a weaker statement (\ref{ub10}). Fortunately, relation (\ref{ub10}) is still sufficient to carry over the proof of (\ref{conv}) in the same way as in \cite{FNP}.

Thus our task is to show that
\bFormula{wws2}
\lim_{n \to \infty} \int_0^T \int_{R^3} \left[
\psi \vre \vue \cdot \Grad \Big( \eta \vc{w}_n \cdot \Grad \Delta^{-1}
[ \xi b(\vre) ] \Big) +
\vre \vc{w}_n \cdot \partial_t \Big(
\psi \eta \Grad \Delta^{-1} [ \xi b(\vre) ] \Big) \right] \dxdt
\eF
\[
= \int_0^T \int_{R^3} \left[
\psi \vr \vu \cdot \Grad \Big( \eta \vc{w} \cdot \Grad \Delta^{-1}
[ \xi \Ov{b(\vr)} ] \Big) +
\vr \vc{w} \cdot \partial_t \Big(
\psi \eta \Grad \Delta^{-1} [ \xi \Ov{b(\vr)} ] \Big) \right] \dxdt,
\]
which, after a straightforward manipulation, reduces to
\bFormula{wws3}
\lim_{n \to \infty}  \int_0^T \int_{R^3}
\psi \eta \left[ \vre \vc{w}_n \cdot \Big( \vue \cdot \Grad \Delta^{-1} \Div [
\xi b(\vre) ] - \Grad \Delta^{-1} \Div [ \xi b(\vre) \vue ] \Big) \right]
 \ \dxdt
\eF
\[
=  \int_0^T \int_{R^3}
\psi \eta \left[ \vr \vc{w} \cdot \Big( \vu \cdot \Grad \Delta^{-1} \Div [
\xi \Ov{b(\vr)} ] - \Grad \Delta^{-1} \Div [ \xi \Ov{b(\vr)} \vu ] \Big) \right]
 \ \dxdt.
\]

Now, with (\ref{ub4}) at hand, we observe that the commutator
\[
\Big( \vue \cdot \Grad \Delta^{-1} \Div [
\xi b(\vre) ] - \Grad \Delta^{-1} \Div [ \xi b(\vre) \vue ] \Big)
\]
is actually \emph{bounded} in $L^2(0,T; W^{1,q} (R^3;R^3))$ for any
$1 < q < 2$. This is easy to see computing ${\rm div}$ and ${\bf curl}$
of this quantity and using continuity of the Riesz operator in $L^q-$spaces, see also the first part of the proof of Theorem 10.28 in \cite{FEINOV}.

On the other hand, as a direct consequence of (\ref{s3}), (\ref{ub5});
\[
\vre \vc{w}_n \to  \vr \vc{w} \ \mbox{in} \ L^2(0,T; W^{-1,q'}(\Omega;R^3));
\]
whence (\ref{wws3}) follows. We have proved (\ref{wws1}).

\subsection{Conclusion}

Having established (\ref{ws1}) we can finish the proof of (\ref{c1}) by the same arguments as in \cite{FNP}. Moreover, the pointwise convergence of the densities may be strengthened to
\bFormula{wsc}
\vre \to \vr \ \mbox{in}\ C([0,T]; L^q(\Omega)) \ \mbox{for any}\
1 \leq q < \gamma.
\eF
In particular, $p(\vr) = \Ov{p(\vr)}$. We have proved Proposition \ref{Ps1}.

\section{Existence of solutions for problems with irregular forces}
\label{e}

Our aim is to use Proposition \ref{Ps1} to show \emph{existence} of
weak solutions to the Navier-Stokes system driven by $\vc{w}$,
\bFormula{e1}
{\rm ess} \sup_{t \in (0,T)} \| \vc{w} (t, \cdot) \|_{W^{1,\infty}(\Omega;R^3)}
\leq C_w.
\eF
To this end, suppose first that $\vc{w}$ is regular, more
specifically,
\bFormula{e2}
\vc{w} \in C^1([0,T]; W^{1,\infty}(\Omega;R^3)),\ \vc{w}(0,\cdot) =
0.
\eF

As shown in \cite[Theorem 1.1]{FNP}, the Navier-Stokes system admits
a weak solution in the sense specified at the beginning of Section
\ref{s}, where the momentum equation (\ref{Fws5}) is replaced by the
integral identity
\bFormula{e3}
\int_0^T \intO{ \Big( \vr \vu \cdot \partial_t \varphi + \vr \vu
\otimes \vu : \Grad \varphi + p(\vr) \Div \varphi \Big)  } \ \dt
\eF
\[
\int_0^T \intO{ \Big( \tn{S}(\Grad \vu ) : \Grad \varphi(x)  - \vr
\partial_t \vc{w} \cdot \varphi \Big)} \ \dt - \intO{ (\vr \vu)_0 \cdot \varphi (0,\cdot) }
\]
for all $\varphi \in \DC([0,T) \times \Omega;R^3)$, while the
corresponding energy inequality reads
\bFormula{e4}
- \int_0^T \intO{ \left( \frac{1}{2} \vr | \vu  |^2 + P(\vr )
\right)} \ \partial_t \psi \ \dt + \int_0^T \intO{ \tn{S} (\Grad \vu
): \Grad \vu } \ \psi \ \dt
\eF
\[
\leq \psi(0) \intO{ \left( \frac{1}{2} \frac{ |(\vr \vu)_0 |^2
}{\vr_0}  + P(\vr_0)  \right) } + \int_0^T \intO{ \vr \partial_t
\vc{w} \cdot \vu } \ \psi \ \dt
\]
for any $\psi \in \DC[0, T)$, $\psi \geq 0$. The result, as stated
in \cite{FNP}, requires certain smoothness of the boundary $\partial
\Omega$ that can be relaxed, see Kuku\v cka \cite{Kuk}. In
particular, the existence can be shown provided $\Omega$ is a
bounded Lipschitz domain.

By virtue of the identity
\[
\int_0^T \intO{ \vr \partial_t \vc{w} \cdot \varphi } \ \dt  =
\int_0^T \intO{ \vr
\partial_t (\vc{w} \cdot \varphi )} \ \dt - \int_0^T \intO{ \vr \vc{w} \cdot
\partial_t \varphi}\ \dt
\]
\[
= - \int_0^T \intO{ \vr \vu \cdot \Grad (\vc{w} \cdot \varphi ) } \
\dt - \int_0^T \intO{ \vr \vc{w} \cdot
\partial_t \varphi}\ \dt
\]
satisfied for any $\varphi \in \DC([0,T) \times \Omega;R^3)$,
relation (\ref{e3}) reduces to (\ref{Fws5}).

Moreover, taking $\varphi = \vc{w}$ as a test function in (\ref{e3})
yields
\bFormula{e5}
\int_0^\tau \intO{ \Big( \vr \vu \cdot \partial_t \vc{w} + \vr \vu
\otimes \vu : \Grad \vc{w} + p(\vr) \Div \vc{w} \Big)  } \ \dt
\eF
\[
= \int_0^\tau \intO{ \Big( \tn{S}(\Grad \vu ) : \Grad \vc{w}  - \vr
\frac{1}{2} \partial_t |\vc{w}|^2 \Big)} \ \dt + \intO{ (\vr \vu)
(\tau, \cdot ) \cdot \vc{w}(\tau, \cdot) }
\]
\[
= - \intO{ \left( \frac{1}{2} \vr |\vc{w}|^2 - \vr \vu \cdot \vc{w}
\right)(\tau, \cdot) } + \int_0^\tau \intO{ \left( \frac{1}{2} \vr
\vu \cdot \Grad |\vc{w}|^2 + \tn{S} (\Grad \vu) : \Grad \vc{w}
\right) } \ \dt \ \mbox{for any}\ \tau \geq 0,
\]
which, together with (\ref{e4}), gives rise to the energy inequality
(\ref{Fws6}).

Having obtained solutions for any regular $\vc{w}$ we may
approximate a general function $\vc{w}$ by a sequence of smooth
functions $\vc{w}_n$ satisfying (\ref{s3}) and use Proposition
\ref{Ps1} to conclude:

\bProposition{e1} Let $\Omega$ be a bounded Lipschitz domain.
Suppose that the pressure $p$ satisfies (\ref{w4}),
\bFormula{e6}
\vc{w} \in L^\infty_{\rm weak}(0,T; W^{1,\infty}_0(\Omega;R^3)),
\eF
and
\bFormula{e7}
\vr_0 \in L^\gamma (\Omega) ,\ \vr_0 \geq 0,\ \intO{ \vr_0 } = M >
0,\ (\vr \vu)_0 \in L^1(\Omega;R^3),\ \intO{ \frac{ |(\vr \vu)_0 |^2
}{\vr_0} } < \infty.
\eF

Then the Navier-Stokes system admits a weak solution in $(0,T)
\times \Omega$ in the sense specified in (\ref{Fws2} - \ref{Fws6}).
\eP

{\bf Remark \ref{e}.1} {\it The symbol $L^\infty_{\rm weak}$ denotes
the space of weakly measurable functions.}

\medskip

{\bf Remark \ref{e}.2} {\it Writing
\[
(\vr \vu)_0 = \sqrt{\vr_0} \frac{ (\vr \vu)_0 }{ \sqrt{\vr_0}}
\]
we can see that (\ref{e7}) implies that
\[
(\vr \vu)_0 \in L^{2\gamma/(\gamma + 1)}(\Omega;R^3).
\]
}

\section{Problems driven by stochastic forces}
\label{ps}

We start by identifying the function space
$\mathcal{W}$ for $\vc{w}$. We suppose that
\begin{enumerate}
\item
$\mathcal{W}$ is a separable complete metric space;
\item
\bFormula{ps1a}
\mathcal{W} \subset L^\infty_{\rm weak} (0,T; W^{1,\infty}_0(\Omega;R^3));
\eF
\item
if
\[
w_n \to w \ \mbox{in}\ \mathcal{W},
\]
then
\bFormula{ps2a}
{\rm ess} \sup_{t \in (0,T)} \| \vc{w}_n (t, \cdot) \|_{W^{1, \infty}_0(\Omega;R^3)}
\leq c , \ \mbox{and}\ w_n \to w \ \mbox{weakly - (*) in} \ L^\infty(0,T; W^{1,\infty}(\Omega)).
\eF

\end{enumerate}

Next, we introduce the spaces of the data. In accordance with
hypotheses of Proposition \ref{Pe1}, we take
\bFormula{ps1}
\mathcal{ID} = \left\{ [\vr_0, (\vr \vu)_0] \Big| \ \vr_0 \in
L^\gamma(\Omega),\ (\vr \vu)_0 \in L^{2\gamma/(\gamma +
1)}(\Omega;R^3), \right.
\eF
\[
\left. \vr_0 \geq 0,\ \intO{\vr_0 } = M > 0,\ \intO{ \frac{1}{2}
\frac{|(\vr \vu)_0 |^2}{\vr_0} + P(\vr_0) } \leq E_0 \right\}.
\]
We note that $\mathcal{ID}$ is a closed convex subset of the
separable Banach space $L^\gamma(\Omega) \times L^{2\gamma/(\gamma +
1)} (\Omega;R^3)$, in particular, it is a Suslin space.

We suppose that the mapping
\[
\omega \in \mathcal{O} \mapsto [\vr_0 (\cdot, \omega), (\vr \vu)_0
(\cdot, \omega) ] \in \mathcal{ID}
\]
is a random variable on $\mathcal{O} = \{ \mathcal{O}, \mathcal{B},
\mu \}$, with a regular probability measure $\mu$. Thus for any $\ep
> 0$, there exists a compact set $\mathcal{K}_\ep \subset
\mathcal{O}$ such that
\[
\mu \left( \mathcal{O} \setminus \mathcal{K}_\ep \right) < \ep,
\]
and the mapping
\[
\omega  \mapsto [\vr_0 (\cdot, \omega), (\vr \vu)_0 (\cdot, \omega)]
\ \mbox{restricted to}\ \mathcal{K}_\ep \ \mbox{is continuous.}
\]

Similarly, we assume that
\[
\omega \in \mathcal{O} \mapsto \vc{w}(\cdot, \omega) \in \mathcal{W}
\]
is a random variable.

Finally, we consider a multi-valued mapping
\[
\mathcal{S}:  \Big\{ [\vr_0, (\vr \vu)_0], \vc{w} \Big\} \in
\mathcal{ID} \times \mathcal{W} \mapsto [\vr, \vu] \in C([0,T];
L^1(\Omega)) \times L^2(0,T;W^{1,2}_0 (\Omega;R^3))
\]
that assigns to the initial data $\vr_0$, $(\vr \vu)_0$ and to the
forcing $\vc{w}$ the associated family of weak solutions to the
Navier-Stokes system, the existence of which is guaranteed by
Proposition \ref{Pe1}. Similarly to the situation examined by Bensoussan and Temam \cite{BenTem}, the
weak solution need not (is not known to be) be unique. As a consequence of Proposition \ref{Ps1},
the values of $\mathcal{S}$ are non-empty closed subsets of the
Banach space $C([0,T]; L^1(\Omega)) \times L^2(0,T;W^{1,2}_0
(\Omega;R^3))$, and $\mathcal{S}$ possesses a closed graph.

Applying the abstract result of Bensoussan and Temam \cite[Theorem
3.1, Lemma 3.1]{BenTem} we conclude that for each $\omega \in
\mathcal{O}$, there exists a weak solution of the Navier-Stokes
system in the sense specified through (\ref{ws1} - \ref{ws5}) such
that the mapping
\[
\omega \mapsto [\vr(\cdot, \omega), \vu(\cdot, \omega)] \in C([0,T];
L^1(\Omega)) \times L^2(0,T; W^{1,2}_0 (\Omega;R^3))
\]
is a random variable.

We have shown the main result of the present paper:
\bTheorem{ps1}
Let $\Omega \subset R^3$ be a bounded Lipschitz domain. Suppose that the pressure $p = p(\vr)$ satisfies hypothesis (\ref{w4}), and that $\mathcal{W}$ is a complete separable metric space satisfying (\ref{ps1a}), (\ref{ps2a}).
Let
\[
\omega \in \mathcal{O} \mapsto \Big\{ [\vr_0(\cdot, \omega), (\vr \vu)_0 (\cdot, \omega)] , \vc{w} \Big\} \in \mathcal{ID} \times \mathcal{W}
\]
be a random variable.

Then there exists a random variable
$\vr (\cdot, \omega) , \ \vu (\cdot, \omega)$,
\[
\vr(\cdot, \omega) \in C([0,T]; L^1(\Omega)), \ \vu (\cdot , \omega) \in L^2(0,T; W^{1,2}_0(\Omega))
\ \mbox{for a.a.}\ \omega \in \mathcal{O},
\]
satisfying the Navier-Stokes system (\ref{ws1} - \ref{ws6}) for a.a. $\omega \in \mathcal{O}$.
 \eT

\section{Applications}
\label{sk}

We conclude the paper by presenting a more specific type of random perturbation that satisfies the general assumptions stated in Section \ref{ps}, which is the so-called L\'evy noise. For the convenience of the reader, the definition and basic properties of this type of random perturbations as well as some basic examples fitting into this framework are recalled. For more details see e.g. the basic monographs \cite{PZ} or \cite{Apple}.

Let $(\mathcal O, \mathcal F, \mathbb P)$ be a probability space and $E = (E,\|\cdot \|_E)$ a separable Banach space.

\begin{Definition}\label{def:X1}
An $E$-valued L\'evy process is a stochastic process $L=(L(t),t\ge 0)$, $L: [0, \infty) \times \mathcal O \rightarrow E$ satisfying:
\begin{itemize}
\item[\rm(i)] $L$ has independent increments, meaning for each $0 \le t_1 < t_2 < \dots t_n$, $n\in\mathbb N$, the $E$-valued random variables $L(t_1)-L(t_0)$, $L(t_2) - L(t_1),\dots$, $L(t_n) - L(t_{n-1})$ are stochastically independent.

\item[\rm(ii)] $L$ has stationary increments, meaning the probability laws of $L(t) - L(s)$, $t>s$, depend only on the difference $t-s$.

\item[\rm (iii)] $L(0) = 0$ and $L$ is stochastically continuous, meaning for each $r>0$ and $s\in\mathbb [0, \infty)$,
\bFormula{Y1}
\lim_{t\rightarrow s} \mathbb P[\|L(t)-L(s)\|_E>r] = 0.
\eF
\end{itemize}
\medskip
\end{Definition}

Note that if the other conditions of Definition~\ref{def:X1} are satisfied, (\ref{Y1}) is equivalent to
\bFormula{Y2}
\lim_{t\rightarrow 0+} \mathbb P[\|L(t)\|_E > r] =0
\eF
for each $r>0$.

Recall that a process is called \emph{c\`adl\`ag} if its  random paths are right-continuous and there exist left limits at each point.

\bTheorem{theo:X2}
Every $E$-valued L\'evy process $L$ has a  c\`adl\`ag modification, that is, there exists a c\`adl\`ag process $\tilde L$ which is a L\'evy process itself (i.e. all random paths of $\tilde L$ are $[0, \infty) \rightarrow E$ c\`adl\`ag) such that
\[
\mathbb P [L(t) = \tilde L(t)] = 1, \ t\in [0, \infty).
\]
\eT

The proof is easily deduced from the classical Kinney's result on c\`adl\`ag modification for Markov processes, cf. \cite{Kin}.

As is usual in probability theory we may identify a stochastic process with its modification, so with no loss of generality we assume in the sequel that $L$ is c\`adl\`ag.
\medskip

\begin{Example}\label{ex:X3}
\textsc{$E$-valued Brownian motion with drift.} Let $L(t) = at + W(t)$, where $a\in E$ is fixed and $W=(W(t),t \in [0, \infty))$ is an $E$-valued centered Gaussian process with the incremental covariance operator $R \in\mathcal L(E^*,E)$, i.e.
\[
\mathbb E \langle W(s),x\rangle_{E,E^*} \langle W(t),y\rangle_{E,E^*} = \langle Rx,y\rangle_{E,E^*}(s \wedge t)
\]
for each $x,y \in E^*$, $s,t\in [0, \infty)$. For example, if $E$ is a Hilbert space, an operator $R\in\mathcal L(E)$ is an incremental covariance of a Wiener process only if $R$ is nonnegative, symmetric and trace class. The process $L$ is then a L\'evy process (with continuous paths). It corresponds to the case when the equation is perturbed by noise white in time.
\end{Example}

\medskip

\begin{Example} \label{ex:X4}
Let $\nu$ be a finite measure on $E\setminus \{0\}$. An $E$-valued compound Poisson process with the L\'evy measure $\nu$ is a c\`adl\`ag L\'evy process satisfying
\[
\mathbb P[L(t)\in A] = e^{-\nu(E)t} \sum^{\infty}_{k=0} \frac{t^k}{k!} \nu^{*k}(A), \quad t \ge 0, \quad \Gamma \in \mathcal B(E) \setminus \{0\}.
\]
An equivalent definition that yields also the construction of the compound Poisson process is this: Let  $(Y_n)_{n\in\mathbb N}$ be a sequence of independent $E$-valued random variables with identical distribution $\frac1{b} \nu$, $b = \nu (E\setminus \{0\})$, and assume that $(\Pi (t), t \ge 0)$ is a standard $\mathbb R_+$-valued Poisson process. Then
\bFormula{Y3}
L(t) = \sum^{\Pi (t)}_{k=0} Z_k
\eF
is a compound Poisson process with the L\'evy measure $\nu$. Clearly, $L$ is a pure jump process with jump size $Z_k$ where $\Pi (t)$ is the number of jumps on the interval $(0,t)$. If $\mathbb E L(t) = t\int_{E\setminus \{0\}} x \nu(\dx) < \infty$, we may define the process $\hat L(t):= L(t) - \mathbb EL(t)$, which is called the compensated compound Poisson process.
\end{Example}
\medskip

The special L\'evy processes introduced in Examples \ref{ex:X3} and \ref{ex:X4} may be used to formulate a useful characterization of a general L\'evy process.

\bTheorem{theo:X5}
{\rm (L\'evy-Khinchine decomposition)}. Every $E$-valued L\'evy process has the representation
\bFormula{Y4}
L(t) = at + W(t) +  \sum^{\infty}_{n=1} L_n (t) + L_0 (t), \ t \ge 0,
\eF
where $a \in E$, $(W(t), t \ge 0)$ is an $E$-valued Brownian motion, $L_n$ for each $n\ge 1$ is a compensated compound Poisson process, $L_0$ is a compound Poisson process, and $W, L_0, L_n$ are stochastically independent. The series in $(\ref{Y4})$ converges $\mathbb P$-almost surely uniformly on each bounded interval of $\mathbb R_+$.
\eT

The processes $L_0,L_n$ may be defined so that their respective L\'evy measures are concentrated on the set $\{ y; \|y\|_E\ge r_0\}$, $\{y; r_{n+1}\le \|y\|_E< r_n\}$, for any arbitrarily chosen sequence $r_n \searrow 0$.

The natural state space for the paths on a time interval $[0,T]$ of $L$ is the linear space of $E$-valued c\`adl\`ag functions $\mathcal D = \mathcal D ([0,T];E)$. It may be equipped with the uniform norm but is not separable in this norm as needed in Section \ref{ps}. A usual way to overcome this difficulty is to consider the so-called Skorokhod metric. Let us briefly recall its definition, for details we may refer to \cite[ Chapter 3]{Bi}. Let $\Lambda$ denote the set of all strictly increasing, continuous mappings of $(0,T]$ onto itself thus if $\lambda \in \Lambda$ then $\lambda(0) = 0$, $\lambda(T) = 1$. Set
\[
d(x,y) = \inf\{\varepsilon >0| \ \mbox{there exists}\  \lambda \in \Lambda ;\lambda(t) - t| < \varepsilon, \ \|x(t) - y(\lambda(t))\|_E < \varepsilon
\ \mbox{for each} \ t\in [0,T]\}.
\]
It may be shown that the space $(\mathcal D,d)$ is separable but not complete, however, it is possible to introduce a metric $d_0$ on $\mathcal D$, equivalent to $d$, such that the space $(\mathcal D,d_0)$ is both separable and complete. Obviously, we have that
\bFormula{ex:Y5}
x_n \rightarrow x \ \mbox{in} \ (\mathcal D, d_0) \Leftrightarrow \ \mbox{there exists}\ \lambda_n \in\Lambda, \ \sup_t (\|x_n(\lambda_n(\cdot))-x(\cdot)\|_E, \ |\lambda_n (t) - t|)\rightarrow 0.
\eF
The conditions imposed on the $\mathcal W$ in Theorem \ref{Tps1} invoke that $E$ should be chosen as a separable Banach space continuously embedded into $W^{1,\infty}_0 (\Omega;\mathbb R^3)$. Thus we may take $E= W^{k,p} \cap W^{1,1}_0 (\Omega,\mathbb R^3)$, where $p\in (1,\infty)$ and $k \ge 1$ is such that
\bFormula{Y6}
k > \frac3{p} +1.
\eF
Thus we report the following result:

\bTheorem{theo:Y6}
Let $\Omega \subset \mathbb R^3$ be a bounded Lipschitz domain. Assume that $p = p(\rho)$ satisfies hypothesis $(\ref{w4})$ and consider the Navier-Stokes system
$(\ref{ws1})$--$(\ref{ws6})$
where $\vc{w}$ is an  $E$-valued L\'evy process, with
\[
E = W^{k,p} \cap W^{1,1}_0 (\Omega,\mathbb R^3), \ k > \frac{3}{p} + 1, \ k \geq 1,\ p > 1.
\]
If the initial datum $(\rho_0,(\rho u)_0)$ is an $\mathcal{ID}$-valued random variable then there exists a random variable $(\rho,u)$ satisfying the Navier-Stokes system.
\eT

The proof follows directly from Theorem \ref{Tps1} if we choose $\mathcal W = \mathcal D = \mathcal D([0,T],E)$. It is easy to check that the conditions (\ref{ps1a}), (\ref{ps2a}) are satisfied.



\def\ocirc#1{\ifmmode\setbox0=\hbox{$#1$}\dimen0=\ht0 \advance\dimen0
  by1pt\rlap{\hbox to\wd0{\hss\raise\dimen0
  \hbox{\hskip.2em$\scriptscriptstyle\circ$}\hss}}#1\else {\accent"17 #1}\fi}

\end{document}